\documentclass[12pt,reqno]{amsart}
\usepackage{amscd,amssymb,amsthm,euscript}

\newtheorem{thm}{Theorem}[section]

\newcounter{ictr}
\newenvironment{ilist}{\begin{list}
                         {(\roman{ictr})}
                         {\usecounter{ictr}
                          \parsep=0pt \itemsep=0pt 
                          \topsep=6pt \partopsep=\parskip
                          \setlength{\leftmargin}{0.5truein}}}
                      {\end{list}}
  
\renewcommand{\epsilon}{{\varepsilon}}
\newcommand{\zbar}{{\overline z}}
\newcommand{\ga}{\Gamma}

\newcommand{\A}{\mathcal{A}}
\newcommand{\B}{\mathcal{B}}

\newcommand{\C}{\mathbb{C}}
\newcommand{\R}{\mathbb{R}}

\renewcommand{\H}{\mathcal{H}}

\newcommand{\x}{\times}
\newcommand{\crg}{C_r^*(\ga)}

\newcommand{\cp}{\rtimes}

\newcommand{\te}{\otimes}

\newtheorem{theorem}{Theorem}[section]
\newtheorem{definition}[theorem]{Definition}

\begin{document}

\title{Exactness and the Novikov Conjecture}
\author{Erik Guentner}
\address{Department of Mathematical Sciences, IUPUI, Indianapolis,
IN 46202-3216}
\email{eguentner@math.iupui.edu}
\author{Jerome Kaminker}
\address{Department of Mathematical Sciences, IUPUI, Indianapolis,
IN 46202-3216} 
\email{kaminker@math.iupui.edu}

\thanks{The first author was supported in part by NSF Grant DMS-9706960}
\thanks{The second author was supported in part by NSF Grant DMS-9706817}

\date{\today} 
\begin{abstract}
  In this note we will study a connection between the conjecture that
  $C^{*}_{r}(\Gamma)$ is exact and the Novikov conjecture for
  $\Gamma$.  The main result states that if the inclusion of the
  reduced $C^*$-algebra $\crg$ of a discrete group $\ga$ into the
  uniform Roe algebra of $\Gamma$, $UC^{*}(\Gamma)$, is a nuclear map
  then $\ga$ is uniformly embeddable in a Hilbert space.  By a result
  of G. Yu, this implies that $\Gamma$ satisfies the Novikov
  conjecture.  Note that the hypothesis is a slight strengthening of
  the usual notion of exactness since a group $\ga$ is exact if and
  only if the inclusion of $\crg$ into $\B(l^2(\ga))$ is nuclear.
\end{abstract}
\maketitle

\section{Introduction}

Let $X$ be a discrete metric space with metric $d$. A function $f$ from
$X$ to a separable Hilbert space $\H$ is a {\it uniform embedding\/} if
there exist non-decreasing proper functions $\rho_\pm : [0,\infty)\to
[0,\infty)$ such that
\begin{eqnarray*}
  \rho_-(d(x,y))\leq \| f(x) - f(y) \| \leq \rho_+(d(x,y)),
   \quad\text{for all $x$, $y\in X$}.
\end{eqnarray*}
  
The Strong Novikov Conjecture states that the assembly map on K-theory,
\begin{equation*}
  \mu : K_{*}(B\Gamma) \to K_{*}(C^{*}_{r}(\Gamma)),
\end{equation*}
is injective.  
Answering a question of Gromov, Yu proved the following theorem
\cite{yu98,skandalis-tu-yu00}. 

\begin{thm}
  Let $\Gamma$ be a finitely presented discrete group.  If $\Gamma$ is
  uniformly embeddable in a Hilbert space, then $\Gamma$ satisfies the
  Strong Novikov conjecture.
\end{thm}
This is currently the weakest general hypothesis implying the Novikov
conjecture.  It is conceivable, however, that there exist groups which
are not uniformly embeddable in a Hilbert space but which nevertheless
satisfy the Novikov conjecture.  At present there are no such examples
known.

From another direction, there is the question of whether all finitely
generated discrete groups are exact,
\cite{kirchberg93,kirchberg-wassermann98}.  Recall that a
discrete group $\Gamma$ is {\it exact\/} if its reduced $C^{*}$-algebra,
$C^{*}_{r}(\Gamma)$, is an exact $C^{*}$-algebra.  That is, given the
exact sequence of $C^{*}$-algebras
\begin{equation*}
  0 \to I \to B \to B/I \to 0
\end{equation*}
the sequence
\begin{equation*}
  0 \to I\te_{min}C^{*}_{r}(\Gamma)  \to B\te_{min}C^{*}_{r}(\Gamma) 
          \to (B/I)\te_{min}C^{*}_{r}(\Gamma) \to 0
\end{equation*}
is also exact.

The Novikov conjecture and Exactness question appear to have little in
common other than that they both involve properties which might be
possessed by all finitely presented groups.  However, there is a link
provided by results of Roe-Higson and Yu \cite{higson-roe98,yu98}.
Combined these results state that if $\Gamma$ acts amenably, in the
topological sense \cite{delaroche-renault98}, on a compact space then it
is uniformly embeddable in a Hilbert space, and hence satisfies the
Novikov conjecture.  On the other hand, it is an easy observation that
this condition also implies that $\Gamma$ is exact.  Thus, the same
hypothesis yields both properties.

We note that Gromov has asserted the existence of finitely presented
groups that are not uniformly embeddable \cite{gromov-s-and-q}.  This
follows from his assertion that there exists a finitely presented group
whose Cayley graph contains a sequence of expanding graphs
\cite{lubotzky-expanders}, together with his observation that, when
viewed as a discrete metric space, a sequence of expanding graphs is not
uniformly embeddable.  On the other hand, it follows simply from a
result of Voiculescu \cite{voiculescu90} that the uniform
algebra of such a metric space is, in general, not exact.  Based on this
it seems likely that Gromov's examples of non-uniformly embeddable
groups will in general fail to be exact.

The purpose of this note is to study the relationship between uniform
embeddability and exactness.  We state the main result, leaving precise
definitions for later in the paper.
We need the uniform Roe algebra, $UC^{*}(\Gamma)$, sometimes called the
``rough'' algebra, introduced by Roe, ~\cite{higson-roe98}.  It is
isomorphic to the reduced cross product $C(\beta\Gamma) \cp_r \Gamma$.
The left regular representation provides an inclusion of $\crg$ into
$\B(l^2(\ga))$, and in fact into $UC^*(\ga)$.  Recall that a group
$\Gamma$ is exact if this inclusion is a nuclear embedding of $\crg$
into $\B(l^2(\ga))$ \cite{wassermann-exactness,kirchberg93}.  We modify
this condition by requiring that the inclusion be a nuclear embedding of
$\crg$ into $UC^*(\ga)$.  Note that most classes of groups which are
know to be exact, including word hyperbolic groups, discrete subgroups
of connected Lie groups, Coxeter groups, etc., actually satisfy this
stronger condition.

\begin{thm}
  Let $\Gamma$ be a finitely presented group.  If the inclusion of the
  reduced $C^*$-algebra $\crg$ into the uniform Roe algebra
  $UC^{*}(\Gamma)$ is a nuclear map, then $\Gamma$ is uniformly
  embeddable in a Hilbert space, and hence satisfies the Novikov
  conjecture.
\end{thm}

It is natural to consider other refinements of exactness which can be
obtained by replacing $UC^{*}(\Gamma)$ by other subalgebras of
$\B(l^{2}(\Gamma))$.  This will be discussed in a future paper.

\section{Approximate units and negative type functions}

In this section we will assemble some of the facts needed for the
results in Section 3.  In particular, we will establish an analog
 of a theorem of Akemann-Walter, ~\cite{ akemann-walter81}. 

  A complex-valued function $f$ on the set $X\times X$ is said to be {\it
  positive definite\/} if, for any $n \geq 1$,
\begin{eqnarray*}
  \sum_{ij} \zbar_i f(x_i,x_j) z_j \geq 0,\quad
       \text{for all $x_1,\dots,x_n\in X$ and $z_1,\dots,z_n\in\C$}.
\end{eqnarray*}
A real-valued function $h$ on $X\times X$ is of {\it 
  negative type\/} if 
\begin{ilist}
  \item $h(x,x)=0$ for all $x\in X$,
  \item $h(x,y)=h(y,x)$ for all $x$, $y\in X$, and
  \item $\sum_{ij} a_i h(x_i,x_j) a_j \leq 0$,
              for all $x_1,\dots,x_n\in X$ and 
         $a_1,\dots,a_n\in\R$ satisfying $\sum_j a_j = 0$.
\end{ilist}

It will be convenient to have the following notation.  If $X$ is a
metric space and $A$ is a subspace, then $C_{0}(X;A)$ will be the set of
functions which tend to zero off of $A$.  That is, $f \in C_{0}(X;A)$ if
for any $\epsilon > 0 $ there is an $R > 0$ such that $|f(x)| <
\epsilon$ if $d(x,A) > R$.

Suppose now that $X$ is a discrete metric space.  Consider the ideal
$C_{0}(X \x X;\Delta) \subseteq l^\infty(X\times X)$ , where $\Delta$ denotes the diagonal of $X \x X$. 
A sequence $f_n\in l^\infty(X\times X)$ satisfies $\|f_nf -f \|\to 0$
for all $f\in C_{0}(X \x X;\Delta)$ if and only if $f_n\to 1$
uniformly on any set of the form $B_{\Delta}(R) = \{(x,y) : d(x,y) <
R \}$.  Finally, we say that a complex-valued function $f$ on $X\times
X$ is {\it metrically proper\/} if it satisfies that for any $C > 0 $
there is an $R > 0$ such that $|f(x,y)| > C$ if $d(x,A) > R$.

The following result is a generalization of ~\cite[Theorem
10]{akemann-walter81} from the case of groups to that of equivalence
relations.

\begin{thm}
\label{unifemb}
  Let $X$ be a discrete metric space.  There exists an approximate unit
  for $C_{0}(X \x X;\Delta)$ consisting of positive definite functions if and only if
  there exists a metrically proper negative type function on
  $X\times X$.
\end{thm}
\begin{proof}
  Let $\phi$ be a metrically proper negative type function on $X \x X$.  By a
  generalization of Schoenberg's Theorem, ~\cite{ berg-forst-ptolcag}, the
  function $ e^{-t\phi(x,y)}$ is positive definite for any
  $t \geq 0$.  Since $\phi$ is metrically proper, one has, for any $t$, $
  e^{-t\phi(x,y)} \in C_{0}(X \x X;\Delta)$.  On the other hand one
  also has $\lim_{t\to 0} \|e^{-t\phi(x,y)}-1\| = 0$ uniformly on
  $B_{\Delta}(R)$ for any $R>0$.  Thus, $\Phi_{t} = e^{-t\phi}$
  provides the approximate unit  for
  $C_{0}(X \x X;\Delta)$ consisting of positive definite functions.   
  
  For the converse, let $u_{\lambda}$ be an approximate unit
  consisting of positive definite functions.  Since $u_{\lambda} \to
  1$ uniformly on $B_{R}(\Delta)$, there exists an $R$ and
  $\lambda_{0}$ so that $u_{\lambda} > 0$ if $d(x,y) < R$ and $\lambda
  > \lambda_{0}$.  One may thus adjust the approximate unit so that
  $u_{\lambda}(x,x) = 1$ for all $x \in X$.  Now, exactly as in
  ~\cite[Theorem 10]{ akemann-walter81} one extracts a sequence
  $u_{\lambda_{i}}$ such that the function $\sum_{n}
  Re(1-u_{\lambda_{n}})2^{n} $ converges to the required metrically
  proper negative type function.
\end{proof}

We next recall the result of Yu, ~\cite{ yu98}, relating metrically
proper negative type functions to uniform embeddings in a Hilbert
space.
\begin{thm}
  The metric space $X$ is uniformly embeddable in a Hilbert space if
  and only if there exists a metrically proper negative type function
  on $X \x X$.
\end{thm}

Combining these two results we obtain
\begin{thm}
\label{auandue}
  The following are equivalent for the countable discrete metric space $X$.
  \begin{ilist}
  \item $X$ is uniformly embeddable in a Hilbert space.
  \item There is a metrically proper negative type function on $X \x
    X$.
  \item There is an approximate unit for $ C_{0}(X \x X;\Delta) $
    consisting of positive definite functions.  
  \end{ilist}
\end{thm}

In Section 3 we will discuss the relation of this to the Haagerup
property for the groupoid $\beta\Gamma \cp \Gamma$.

\section{Exactness}
\label{sec:exactness}

In this section we restrict $X$ to be a finitely presented group with a
length function determined by a finite, symmetric set of generators.
The length function, $l$, determines a right invariant metric via
$d(s,t) = l(st^{-1})$.  The quasi-isometry type of $(\Gamma,d)$ is
independent of the choice of generators.  We next recall the definition
of the uniform Roe algebra associated to $(\Gamma,d)$.

Consider the set of $A : \Gamma \x \Gamma \to \C$ satisfying
\begin{ilist}
\item there exists $M > 0$ such that $|A(s,t)| \leq M$, for all $s$,
  $t\in\ga$ 
\item there exists $R > 0$ such that $A(s,t) = 0$ if $d(s,t) > R$
\end{ilist}
Each such $A$ defines a bounded operator on $l^{2}(\Gamma)$ via the
usual formula for matrix multiplication:
\begin{equation*}
  A\xi(s) = \sum_{r \in \Gamma} A(s,r)\xi(r),
       \quad\text{for $\xi\in l^2(\ga)$}.
\end{equation*}
These will be referred to as {\it finite width operators}. The
collection of finite width operators is a $*$-subalgebra of
$\B(l^2(\ga))$.  The {\it uniform Roe algebra} of $\Gamma$, denoted
$UC^{*}(\Gamma)$, is the closure of the $*$-algebra of finite width
operators.  It is a $C^{*}$-algebra The quasi-isometry class of
$(\Gamma,d)$ determines $UC^*(\ga)$, which is therefore independent of
the choice of generators.

Every $t\in\ga$ acts on $l^2(\ga)$ by the left regular representation.
The action of $t \in \Gamma$ is represented by the matrix $A$ defined by
$A(s,r) = 1$ if and only if $s = tr$.  Clearly, $t\in\ga$ acts as a
finite width operator.  Thus, $\C [ \Gamma ] \subseteq UC^{*}(\Gamma)$,
and we have
\begin{eqnarray*}
  \crg \subseteq UC^{*}(\Gamma) \subseteq \B(l^2(\ga)).
\end{eqnarray*}

Recall that if a unital *-homomorphism, $T: \A \to \B$, between unital
$C^{*}$-algebras is {\it nuclear} then there is a net $T_{\lambda} : \A
\to \B$ of finite rank, unital, completely positive linear maps such
that $\lim_{\lambda} \| T_{\lambda}(x) - T(x)\| = 0$ for all $x \in \A$.
It was shown by Kirchberg \cite{wassermann-exactness,kirchberg95} that a
unital $C^{*}$-algebra $\A$ is exact if and only if every
non-degenerate, faithful representation of $\A$ on a Hilbert space $\H$
provides a nuclear embedding of $\A$ into $B(\H)$.  In particular, a
discrete group $\ga$ is exact if and only if the inclusion of $\crg$
into $\B(l^2(\ga))$ given by the left regular representation is a
nuclear embedding.  The main theorem of this section states that if one
restricts the range of the nuclear embedding a little bit, then this
strengthened form of exactness implies the uniform embeddability of
$\Gamma$.

\begin{thm}
\label{thm:main}
  Let $\ga$ be a finitely generated discrete group.  If the inclusion
  $\crg\subset UC^{*}(\Gamma)$ is a nuclear map then $\ga$ is uniformly
  embeddable in a Hilbert space (and hence satisfies the Novikov conjecture).
\end{thm}
\begin{proof}
  By Theorem~\ref{auandue} it is sufficient to produce an approximate
  unit for $C_{0}(X \x X;\Delta) $ consisting of positive definite
  functions.  This will be obtained using nuclearity of the inclusion.
  
  There is a general procedure to associate to a linear map
  $T:\crg\to\B(l^2(\ga))$ a function $u: \Gamma \x \Gamma \to \C$
  given by the formula
\begin{eqnarray*}
  u(s,t) = 
      \langle \delta_s, T(st^{-1})\delta_t \rangle,
\end{eqnarray*}
where $\delta_{t}$ denotes the characteristic function of the element $t
\in \Gamma$.  Note that if $T$ is bounded then $u\in l^\infty \ga\x\ga$.  
The correspondence
\begin{equation*}
  \{T : C^{*}_{r}(\Gamma) \to \B(l^{2}(\Gamma)) \} \longmapsto 
          \{u : \Gamma \x \Gamma \to \C \}
\end{equation*}
has the following properties:
\begin{ilist}
  \item if $T$ is unital and completely positive then $u$ is positive
  definite, and
  \item if $T:\crg\to UC^{*}(\Gamma)$ has finite rank then $u\in
    C_{0}(\ga\times\ga;\Delta)$.
\end{ilist}
Further if $T_{\lambda} : \crg\to\B(l^2(\ga))$ is a net of bounded
linear maps with associated functions $u_{\lambda}$ then
\begin{ilist}
\setcounter{ictr}{2}
  \item if $\|T_{\lambda}(x)-x\| \to 0$, for all $x\in\crg$, then
    $u_{\lambda}\to 1$ uniformly on $B_{R}(\Delta)$ for all $R$. 
\end{ilist}

We verify these properties below, but for now observe that together they
imply the desired result.  Assuming nuclearity of the inclusion of
$\crg$ into $UC^*(\ga)$ we obtain unital completely positive maps
$T_\lambda: \crg\to UC^*(\ga)$ as above.  It follows immediately from
the properties above that the associated functions $u_\lambda\in
l^\infty(\ga\x \ga)$ form the desired approximate unit..

We now turn to the verification of (i)--(iii), beginning with (i).  Let
$s_1,\dots,s_n\in\ga$ and $z_1,\dots,z_n\in\C$.  Define an element of
$\H=\oplus l^2(\ga)$ by $\xi=(z_1\delta_{s_1},\dots,z_n\delta_{s_n})$
and an operator on $\oplus l^2(\ga)$ by the $n\times n$ matrix $A =
[A_{ij}]\in M_n(\B(l^2(\ga)))$ where
\begin{eqnarray*}
    A_{ij} =  T(s_js_i^{-1}).
\end{eqnarray*}
A direct calculation shows
\begin{eqnarray*}
  \sum_{i,j} \zbar_i u(s_i,s_j)z_j &=& \sum_{i,j} \zbar_i 
  \langle \delta_{s_j}, T(s_js_i^{-1})\delta_{s_i} \rangle  
     =  \langle \xi, A\xi \rangle_\H.
\end{eqnarray*}
Thus, it suffices to show that $A$ is positive operator on $\H$.
However, since $T$ is completely positive, this will follow from the
fact that the $n\x n$ matrix $B =[B_{ij}] \in M_n(\B(l^2(\ga)))$, where
$ B_{ij} = s_js_i^{-1}$ defines a positive operator on $\H$.  This is
equivalent to the assertion that 
\begin{eqnarray*}
  \sum_{ij} \langle (s_js_i^{-1})f_i ,f_j \rangle
        = \| (s_1^{-1}f_1,\dots,s_N^{-1}f_n) \|^2 \geq 0 
\end{eqnarray*}
for all $f_1,\dots,f_n\in l^2(\ga)$, which is straightforward.  

We now prove (ii).  Since $T$ has finite rank there exist finitely many
$f_{i}\in \crg^*$ and $S_{i}\in UC^{*}(\Gamma)$ such that $T = \sum
f_{i}S_{i}$.  Since $u$ depends (conjugate) linearly on $T$ it is
sufficient to consider the rank one case where $T(s) = f(s)S$. In this
case,
\begin{eqnarray*}
  | u(s,t) | = | \langle \delta_s,T(st^{-1})\delta_t \rangle|
        = | f(st^{-1})| |\langle \delta_s,S\delta_t\rangle|
        \leq \|f\|_{\crg^*} |\langle \delta_s,S\delta_t\rangle|,
\end{eqnarray*}
and it suffices show that for all $\epsilon>0$ there exists $R>0$
such that $|\langle \delta_s,S\delta_t\rangle| <\epsilon$ provided
$d(s,t)\geq R$.

At this point the requirement that $S \in UC^{*}(\Gamma)$ is needed.
Let $S'$ be a finite width operator such that $\| S-S' \| < \epsilon$.
Then we have
\begin{eqnarray*}
  | \langle S\delta_t,\delta_s \rangle | \leq
       \| S-S' \| + | \langle S'\delta_t,\delta_s \rangle |, 
\end{eqnarray*}
and for large enough $R$, $d(s,t) > R$ forces the last term to be zero.
The result follows.

We conclude the proof by verifying (iii).    
Consider   
\begin{eqnarray*}
  u(s,t)-1 &=& \langle \delta_s,T(st^{-1})\delta_t \rangle 
                       - \langle \delta_s,\delta_s \rangle \\
        &=& \langle \delta_s,T(st^{-1})\delta_t-\delta_s \rangle \\
        &=& \langle \delta_s,(T(st^{-1}) - st^{-1})\delta_t \rangle.
\end{eqnarray*}
Thus, if we have a family $T_{\lambda}$, it follows that
\begin{eqnarray*}
  | u_{\lambda}(s,t)-1 | \leq \| T_{\lambda}(st^{-1}) - st^{-1} \|.
\end{eqnarray*}
To verify the uniform convergence on sets of the form $B_{R}(\Delta)$
note that $d(s,t) < R$ implies that $st^{-1}$ lies in a bounded subset
of $\Gamma$, hence only a finite number of such products are
possible.  Thus, by taking $\lambda$ sufficiently large the right
side can be made as small as necessary.  This completes the proof.
\end{proof}

\section{Approximate units and the Haagerup property}

The results of Section 2 can be used to directly relate the existence
of an approximate unit of positive definite functions to the Haagerup
property for the transformation groupoid $\beta\Gamma \cp \Gamma$.
Here $\beta\Gamma$ is the Stone-Cech compactification of $\Gamma$ and
$\Gamma$ acts on $\beta\Gamma$ on the right by extending right
translation. 

This requires extending the notion of positive definite and negative
type functions to groupoids.  This has been done by
Tu~\cite[Section~3.3]{tu99} in the following way.  We specialize to the
case of a transformation groupoid $X\rtimes\Gamma$, defined as above,
where $X$ is a compact space on which $\Gamma$ acts on the right.  A
complex-valued function $\varphi$ on $X\rtimes\Gamma$ is {\it positive
  definite\/} if
\begin{eqnarray*}
  \sum_{ij} \zbar_i \varphi(x\cdot s_i,s_i^{-1}s_j) z_j \geq 0,\quad
      \text{for all $x\in X$, $s_1,\dots,s_n\in\Gamma$ 
             and $z_i,\dots,z_n\in\C$}.
\end{eqnarray*}
A real-valued function $\psi$ on $X\rtimes\Gamma$ is of {\it negative
  type\/} if
\begin{ilist}
  \item $\psi(x,e) = 0$ for all $x \in X$,
  \item $\psi(x\cdot s,s^{-1}t) = \psi(x\cdot t,t^{-1}s)$, for all
    $x\in X$, $s$, $t\in\Gamma$,
  \item $\sum_{ij} a_i \psi(x\cdot s_i,s_i^{-1}s_j)a_j \leq 0$, for
    all $x\in X$, $s_1,\dots,s_n\in\Gamma$ and $a_1,\dots,a_n\in\R$
    satisfying $\sum_j a_j = 0$.
\end{ilist}

\begin{definition}
  The transformation groupoid $X\rtimes\Gamma$ has the {\em Haagerup
  property} if there exists a proper, negative type function $\psi :
  X\rtimes\Gamma \to \R$.
\end{definition}

If $X\rtimes\Gamma$ has the Haagerup property, then it admits a proper
affine action on a field of Hilbert spaces, ~\cite{ tu99}.  This
latter property, in the case of groups, is call {\em a-T-menability}.
We may now state the main result of this section.

\begin{thm}
Let $\Gamma$ be a discrete group.  The following are equivalent:
\begin{ilist}
  \item $\Gamma$ is uniformly embeddable in a Hilbert space.
  \item The groupoid $\beta\Gamma\rtimes\Gamma$ has the Haagerup property. 
  \item $C_0(\beta\Gamma\rtimes\Gamma)$ admits an approximate unit of positive
    definite functions.
\end{ilist}  
\end{thm}
\begin{proof}

  There is an equivalence of groupoids, $\alpha : \Gamma\rtimes\Gamma
  \to \Gamma\times\Gamma$ given by $\alpha(s,t) = (s, st)$.  Here
  $\Gamma\times\Gamma$ is the trivial groupoid.  The inverse of
  $\alpha$ is $\beta(s,t) = (s,s^{-1}t)$.  These maps  define
  correspondences, $\alpha^{*}: l^\infty(\Gamma\times\Gamma)
  \leftrightarrow C_b(\beta\Gamma\rtimes\Gamma): \beta^{*}$ between functions on $\Gamma\times\Gamma$ and
  $\beta\Gamma\rtimes\Gamma$ via $\alpha^{*}(f)(s,t) = f(s,st)$ and $\beta^{*}( g)(s,t) = g(s,s^{-1}t)$.
Note that $\alpha^{*}(f)$, initially defined on
$\Gamma\times\Gamma\subset\beta\Gamma\rtimes\Gamma$, extends by continuity to
$\beta\Gamma\rtimes\Gamma$ since $f\alpha (\cdot,t)$ is bounded
for each fixed $t \in \Gamma$.

It is easy to check  that $\alpha^{*}$ and $\beta^{*}$ are inverses and
provide a bijection between $l^\infty(\Gamma\times\Gamma)$ and
$C_b(\beta\Gamma\rtimes\Gamma)$.  These maps have the following
properties which are direct consequences of the definitions.
\begin{ilist}
  \item A function $f \in l^\infty(\Gamma\times\Gamma)$ is
    metrically proper  if and only if $\alpha^{*}(f)$ is a 
    proper function on $\beta\Gamma\rtimes\Gamma$.
  \item The map $\alpha^{*}$ takes the ideal
    $C_{0}(\Gamma\times\Gamma;\Delta)$ to the ideal
    $C_0(\beta\Gamma\rtimes\Gamma)$.
  \item A net $\{ u_{\lambda}\}$ is an approximate unit for
    $C_{0}(\Gamma\times\Gamma;\Delta)$ if and only if
    $\alpha^{*}(f)(u_{\lambda})$ is an approximate unit for
    $C_0(\beta\Gamma\rtimes\Gamma)$.  
\end{ilist}

It remains to note that $\alpha^{*}$ preserves positive definite and
negative type functions.  This also follows in a straightforward way
from the above formulas.

Now the result follows from Theorem ~\ref{auandue}.
\end{proof}

\section{Remarks}

A finitely generated discrete group $\ga$ is {\it strongly exact\/} if
the inclusion of $\crg$ into $UC^*(\ga)$ given by the left regular
representation is a nuclear map, that is, if $\ga$ satisfies the
hypethesis of Theorem~\ref{thm:main}.

\begin{enumerate}
\item If a discrete group $\Gamma$ has the property that there is a
  nuclear embedding of $C^{*}_{r}(\Gamma)$ into $B(l^{2}(\Gamma))$ then
  the inclusion given by the left regular representation is also
  nuclear.  It is possible that this inclusion is also a nuclear map
  into $UC^{*}(\Gamma)$.  In other words, it is possible that every
  exact group is strongly exact.  If this is indeed the case then one
  would deduce that an exact group satisfies the Novikov conjecture.
\item One may consider algebras $A$ satisfying
  \begin{equation*}
    C^{*}_{r}(\Gamma) \subseteq A \subseteq UC^{*}(\Gamma)
  \end{equation*}
and impose the requirement that the inclusion of $\crg$ into $A$ be a
nuclear map; if $A=UC^*(\ga)$ then $\ga$ is strongly exact, whereas if
$A=\crg$ then $\crg$ is nuclear.  One obtains a family of conditions
interpolating between strong exactness and nuclearity.  In the case
$A=\crg$ the procedure employed above for constructing a proper negative
type function on $\Gamma \x \Gamma$ actually yields an invariant one,
which descends to $\ga$ showing that 
$\Gamma$ has the Haagerup property.  This gives an alternate account of
the result of Beki, Cherix and Valette~\cite{bekka-cherix-valette95}.

\end{enumerate}


\begin{thebibliography}{10}

\bibitem{akemann-walter81}
C.~Akemann and M.~Walter, \emph{Unbounded negative definite functions},
  Canadian J.~Math. \textbf{4} (1981), 862--871.

\bibitem{delaroche-renault98}
C.~Anantharaman-Delaroche and J.~Renault, \emph{Amenable groupoids}, Preprint,
  1998.

\bibitem{bekka-cherix-valette95}
M.~E. Bekka, P.~A. Cherix, and A.~Valette, \emph{Proper affine isometric
  actions of amenable groups}, {N}ovikov Conjectures, Index Theorems and
  Rigidity (S.~Ferry, A.~Ranicki, and J.~Rosenberg, eds.), vol.~2, London
  Mathematical Society Lecture Notes, no. 226, 227, Cambridge University Press,
  1995, pp.~1--4.


\bibitem{berg-forst-ptolcag}
C.~Berg and G.~Forst, \emph{Potential theory on locally compact {A}belian
  groups}, Ergebnisse der Mathematik und ihrer Grenzgebiete, vol.~87, Springer,
  New York, 1975.

\bibitem{gromov-s-and-q}
M.~Gromov, \emph{Spaces and questions}, Unpublished manuscript, 1999.

\bibitem{higson-roe98}
N.~Higson and J.~Roe, \emph{Amenable actions and the {N}ovikov conjecture},
  Preprint, 1998.

\bibitem{kirchberg93}
E.~Kirchberg, \emph{On non-semisplit extensions, tensor products and exactness
  of group {$C^*$}-algebras}, Invent.~Math. \textbf{112} (1993), 449--489.

\bibitem{kirchberg95}
\bysame, \emph{On subalgebras of the {CAR}-algebra}, Jour.~Funct.~Anal.
  \textbf{129} (1995), 35--63.

\bibitem{kirchberg-wassermann98}
E.~Kirchberg and S.~Wassermann, \emph{Exact groups and continuous bundles of
  {$C^*$}-algebras}, Preprint, 1998.

\bibitem{lubotzky-expanders}
A.~Lubutzky, \emph{Discrete groups, expanding graphs and invariant measures},
  Progress in Mathematics, vol. 125, Birkh\"auser, Boston, 1994.

\bibitem{skandalis-tu-yu00}
G.~Skandalis, J.~L. Tu, and G.~Yu, \emph{Coarse {B}aum-{C}onnes conjecture and
  groupoids}, Preprint.

\bibitem{tu99}
J.~L. Tu, \emph{La conjecture de {B}aum-{C}onnes pur les feuilletages
  moyennables}, {$K$}-Theory \textbf{17} (1999), 215--264.

\bibitem{voiculescu90}
D.~Voiculescu, \emph{Property {T} and approximation of operators}, Bull.~London
  Math.~Soc. \textbf{22} (1990), 25--30.

\bibitem{wassermann-exactness}
S.~Wassermann, \emph{Exact {$C\sp *$}-algebras and related topics}, Lecture
  Note Series, vol.~19, Seoul National University, Seoul, 1994.

\bibitem{yu98}
G.~Yu, \emph{The {C}oarse {B}aum-{C}onnes conjecture for spaces which admit a
  uniform embedding into {H}ilbert space}, To appear in Invent.~Math., 1998.

\end{thebibliography}

\providecommand{\bysame}{\leavevmode\hbox to3em{\hrulefill}\thinspace}

\end{document}